\documentclass[a4paper]{article}
\usepackage{amsmath,amssymb,amsthm,amscd,tabularx,rotating}
\usepackage{diagrams}
\input cyracc.def 
\font\tencyr=wncyr10

\def\russe{\tencyr\cyracc}
\def\Sha{\text{\russe{Sh}}}

\newcommand{\QQ}{\mathbb Q}

\newcommand{\ZZ}{\mathbb Z}

\DeclareMathOperator{\Sel}{Sel}
\DeclareMathOperator{\coker}{coker}
\DeclareMathOperator{\Hom}{Hom}
\DeclareMathOperator{\rk}{rk}
\DeclareMathOperator{\crk}{cork}
\DeclareMathOperator{\Gal}{Gal}
\DeclareMathOperator{\ord}{ord}

\newcommand{\tE}{\tilde{E}_{v,p^\infty}}
\newcommand{\Ep}{E_{p^\infty}}
\renewcommand{\SE}{\Sel(E/F_\infty)}
\newcommand{\Fiw}{F_{\infty,w}}
\newcommand{\tEE}{\tilde{E}_{v',p^\infty}}
\newcommand{\kiw}{k_{\infty,w}}
\newcommand{\Liw}{L_{\infty,w}}
\newcommand{\Kiw}{K_{\infty,w}}

\newtheoremstyle{mythm}{10pt}{10pt}{\it}{}{\bf}{.}{\newline}{}
\newtheorem{thm}{Theorem}[section]

\theoremstyle{plain}
\newtheorem{prop}[thm]{Proposition}
\newtheorem{lem}[thm]{Lemma}
\newtheorem{cor}[thm]{Corollary}
\newtheorem{con}[thm]{Conjecture}
\newtheorem{rem}[thm]{Remark}

\title{Selmer groups over $p$-adic Lie Extensions I}
\author{Sarah Livia Zerbes}
\date{\today}

\begin{document}  
\maketitle


\section{Introduction}

 \subsection{Statement of the main result}
  Let $F$ be a finite extension of $\QQ$, and let $p$ be an odd prime number. Let $F_\infty$ denote a $p$-adic Lie extension of $F$, i.e. a Galois extension
  of $F$ whose Galois group $\Sigma$ is a $p$-adic Lie group of positive dimension. Perhaps the simplest example of such a $p$-adic Lie extension is the
  cyclotomic $\ZZ_p$-extension of $F$, which we denote by $F^{cyc}$, and which is, by definition, the unique $\ZZ_p$-extension of $F$ contained in the field
  obtained by adjoining all $p$-power roots of unity to $F$. Let us assume that our general $p$-adic Lie extension satisfies the following hypotheses:-
  
  \vspace{1ex}
  
  \noindent {\bf Hypothesis} $\bf{I}$
  \begin{enumerate}
   \item $F_\infty$ is unramified outside a finite set of primes of $F$,
   \item $F_\infty$ contains the cyclotomic $\ZZ_p$-extension $F^{cyc}$ of $F$,
   \item $\Sigma$ has no element of order $p$.
  \end{enumerate} 
  \vspace{1ex}

  A basic result of Lazard~\cite{lazard} and Serre~\cite{serre4} shows that, by virtue of (iii),
   $\Sigma$ has finite $p$-cohomological dimension, which is equal to its dimension
  as a $p$-adic Lie group. Now let $W$ be a discrete $p$-primary module. We say that $W$ has finite $\Sigma$-Euler characteristic if the
  $H^i(\Sigma,W)$ are finite for all $i\geq 0$. When $W$ has finite $\Sigma$-Euler characteristic, we define, as usual, its Euler
  characteristic $\chi(\Sigma,W)$ by
  \begin{equation}
   \chi(\Sigma,W)=\prod_{i\geq 0}(\# H^i(\Sigma,W))^{(-1)^i}.
  \end{equation}
  This paper will be concerned with the study of the Euler characteristics of Selmer groups over $F_\infty$ of elliptic curves, which are defined
  over $F$. We recall that, for any algebraic extension $K$ of $F$, the Selmer group of $E$ over $K$, denoted by $\Sel(E/K)$, is defined by
  \begin{equation}
   \Sel(E/K)=\ker(H^1(K,\Ep) \rightarrow \prod_w H^1(K_w,E)),
  \end{equation}
  where $\Ep$ denotes the group of all $p$-power division points on $E$, and the groups are the usual Galois cohomology groups with the one minor caveat
  that, if $K$ is an infinite extension of $F$, then $K_w$ denotes as usual the union of completions at $w$ of all finite extensions of $F$ contained in $K$.
  If $K$ is Galois over $F$, the natural action of $\Gal(K/F)$ on $H^1(K,\Ep)$ induces an action of $\Gal(K/F)$ on $\Sel(E/K)$. 
  
  What can we say about the $\Sigma$-Euler characteristic of $\SE$? Under fairly general circumstances, which we now explain, it seems likely that $\SE$
  has finite $\Sigma$-Euler characteristic which is very closely related to the $\Gamma$-Euler characteristic of $\Sel(E/F^{cyc})$, where
  $\Gamma=\Gal(F^{cyc}/F)$. Specifically, let us assume, in addition to Hypothesis I, that 
  
  \vspace{1ex}
  \noindent {\bf Hypothesis} $\bf{II}$ 
  \begin{enumerate}
  \item $\Sel(E/F)$ is finite,
  \item $E$ has good ordinary reduction at all primes $v$ of $F$ dividing $p$.
  \end{enumerate}
  \vspace{1ex}
    
  Under these assumptions, Schneider and Perrin-Riou have shown that there is a simple explicit formula of Birch-Swinnerton-Dyer type for $\chi(\Gamma,\Sel(E/F^{cyc}))$
  (see~\cite{coatessujatha1} and formula \eqref{cyclotomic} below). It does not seem unreasonable to hope that, more generally, $\SE$ has finite $\Sigma$-Euler
  characteristic, which is of the form
  \begin{equation}\label{formulaEC}
   \chi(\Sigma,\SE)=\chi(\Gamma,\Sel(E/F^{cyc}))\times\mid\prod_{w\in\frak{M}}L_v(E,1)\mid_p.
  \end{equation}
  Here, $\frak{M}$ denotes a subset of primes of $F$ which do not divide $p$ but which ramifiy in $F_\infty$, $L_v(E,s)$ denotes the Euler factor at $v$ of
  the complex $L$-function of $E$ over $F$ and $\mid$ $\mid_p$ is the $p$-adic valuation of $\QQ$, which is normalized so that $\mid p\mid_p=p^{-1}$. 
  In order to state the formula for $\chi(\Gamma,\Sel(E/F_{cyc}))$, we have to introduce the following notation:
  Let $v$ be a finite prime of $F$, $F_v$ the completion of $F$ at $v$ and $k_v$ the residue field of $v$. Put $c_v=[E(F_v):E_0(F_v)]$, 
  where $E_0(F_v)$ is the subgroup of $E(F_v)$ consisting of the points of non-singular reduction modulo $v$. Write $\tilde{E}_v$ for the reduction of $E$
  modulo $v$. Let $\Sha(E/F)$ be the Tate-Shafarevich group of $E$ over $F$, and $\Sha(E/F)(p)$ its $p$-primary subgroup. Then, under the assumption that $E$ has
  good ordinary reduction at all primes $v$ of $F$ dividing $p$, and that $\Sel(E/F)$ is finite, the classical formula for $\chi(\Gamma,\Sel(E/F^{cyc}))$
  asserts that
  \begin{equation}\label{cyclotomic}
   \chi(\Gamma,\Sel(E/F^{cyc}))=\rho_p(E/F),
  \end{equation}
  where
  \begin{equation}
   \rho_p(E/F)=\frac{\#\Sha(E/F)(p)}{\#(E(F)(p))^2}\times\left\vert\prod_vc_v\times\prod_{v\mid
   p}\#(\tilde{E}_v(k_v))^2\right\vert_p^{-1}
  \end{equation}
  
  \vspace{1ex}
  In support of the general question above, the following special cases of the $\Sigma$-Euler characteristic of $\SE$ have been studied already (always assuming Hypotheses I and
  II):-
  \vspace{1ex}
  \vspace{1ex}
  \begin{enumerate}
   \item $F_\infty=F(\Ep)$, where $E$ is an elliptic curve defined over $F$ with complex multiplication by Perrin-Riou ~\cite{perrin-riou}
   \item $F_\infty=F(\Ep)$, where $E$ is an elliptic curve defined over $F$ without complex multiplication (the so-called $GL_2$-extension) by Coates and
    Howson~\cite{coates2}  (for an alternative approach see~\cite{coatesschneidersujatha})
   \item $F_\infty=F(\mu_{p^\infty},\sqrt[p^\infty]{\alpha})$, where $\alpha \in F^*\backslash\mu$ (the so-called false Tate curve extension) by Hachimori and
   Venjakob~\cite{hachimorivenjakob}
  \end{enumerate}
  \vspace{1ex}
  
  In this paper, we will prove~\eqref{formulaEC} in the case when $F_\infty$
  is obtained by adjoining to $F$ the $p$-division points of any abelian variety defined over $F$. To be precise, we prove the following result:-
  
  \begin{thm}\label{EC}
   Let $A$ be an abelian variety defined over a finite extension $F$ of $\QQ$, let $p$ be a prime $\geq 5$, let
   $F_\infty = F(A_{p^{\infty}})$ and $\Sigma=\Gal(F_\infty/F)$. 
   Let $E$ be an elliptic curve defined over $F$. Suppose that 
   (i) $\Sigma$ does not have any element of order $p$,
   (ii) both $A$ and $E$ have good ordinary reduction at all primes 
   of $F$ dividing $p$, (iii) $\Sel(E/F)$ is finite, and (iv) the Pontryagin dual ${\cal C}(E/F_\infty)$ of $\SE$ is $\Lambda(\Sigma)$-torsion.
   Then $\SE$ has finite $\Sigma$-Euler characteristic, and
   \begin{equation}
    \chi(\Sigma,\SE)=\chi(\Gamma,\Sel(E/F^{cyc})) \times\mid\prod_{v \in \frak{M}} L_v(E,1)\mid _p,
   \end{equation} 
   where $\frak{M}$ is the finite set of primes of $F$ where $A$ has bad, but not potentially good reduction.
  \end{thm}
  
  Here, $\Sigma$ is a closed subgroup of $GL_n(\ZZ_p)$, where $n=2\dim(A)$, and thus is a compact $p$-adic Lie group. By the Weil pairing, the
  cyclotomic $\ZZ_p$-extension is contained in $F_\infty$. Finally, we will show in Section~\ref{local} that a prime $v$ of $F$ ramifies in $F_\infty$ if
  and only if either $v$ divides $p$ or $A$ has bad and not potentially good reduction at $v$, so $F_\infty$ does indeed satisfy Hypothesis I. However, in
  contrast to the cases (i)-(iii) above, no explicit description of $\Sigma$ is known in general, and it is shown in~\cite{pink} that the dimension of
  $\Sigma$ can be arbitrarily large. Note that condition (iii) of Hypothesis I is automatically staisfied if $p>2\dim(A)+1$.
  
  The theorem gives further evidence to the conjecture that under the circumstances of Hypothesis II it should be possible to prove an Euler
  characteristic formula for $\SE$ of the form~\eqref{formulaEC} over a general $p$-adic Lie extension satisfying Hypothesis I. Looking closely at the
  proof of the theorem reveals that most of the arguments work in complete generality, and that the structure of $F_\infty$ as the division field of
  points of finite order on an abelian variety is used only in the proofs of Propositions~\ref{locEC} and~\ref{globEC}. Hence a general proof of the
  conjecture under the assumption of Hypothesis II reduces to showing that the following conditions are satisfied:-\\
  
  \vspace{1ex}
  \noindent {\bf Condition I} 
  $\Ep(F_\infty)$ has finite $\Sigma$-Euler characteristic, and
  \begin{equation}
   \chi(\Sigma,\Ep(F_\infty))=1;
  \end{equation}
  
  \noindent {\bf Condition II} 
  For any prime $v$ of $F$ dividing $p$, $\tE(\kiw)$ has finite $\Sigma_w$-Euler characteristic, and 
  \begin{equation}
   \chi(\Sigma_w,\tE(\kiw))=1.
  \end{equation}
  Even though the highly technical proof of Theorem~\ref{EC} reduces to proving these comparatively simple statements, it seems very difficult to attack
  these problems in general, since the proof would require a detailed analysis of the structures of the global and local Galois groups $\Sigma$ and
  $\Sigma_w$. In particular, it seems fairly unlikely that the arguments which are used in the proofs of Propositions~\ref{locEC} and~\ref{globEC}
  extend to larger classes of $p$-adic Lie extensions in any obvious way, since they make crucial use of the fact that $F_\infty$ is the $p$-division
  field of an abelian variety. However, it is remarkable that no counterexamples to the above statements are known.

  \paragraph{Notation.} The following notation is used frequently in this paper:
  For a field $L$, $\bar{L}$ denotes its separable closure. If $L$ is a local field and $M$ is a $\Gal(\bar{L}/L)$-module, then we write $H^i(L,M)$ instead
  of $H^i($Gal$(\bar{L}/L),M)$. Throughout, $S$ will denote a finite set of non-archimedean primes of $F$. We write $F_S$ for the maximal extension of $F$
  unramified outside $S$ and the infinite primes, and for any algebraic extension $H$ of $F$ contained in $F^S$, we let $G_S(H)=$Gal$(F^S/H)$. If $M$ is an abelian group, we write $M(p)$
  for its $p$-primary subgroup. When $M$ is a discrete $p$-primary Abelian group or a compact pro-$p$ group, we define its Pontryagin dual 
  \begin{equation}
   \hat{M}=\Hom(M,\QQ_p/\ZZ_p).
  \end{equation}

  \paragraph{Acknowledgements.} I would like to thank John Coates for his advice and encouragement. Also, thanks go to 
  Konstantin Ardakov, Vladimir Dokchitser, Susan Howson, Burt Totaro and Otmar Venjakob for their interest and for many helpful discussions and to Christian
  Wuthrich for his astonishing patience when teaching me \LaTeX. 
  
 \subsection{Strategy}

 Let $S$ be a finite set of primes of $F$ containing all the primes dividing $p$, all the primes where $E$ has bad reduction and all the primes where 
 $A$ does not have potentially good reduction. Let $L$ be a finite extension of $F$. For each finite place $v$ of $F$, define 

 \begin{equation}
  J_v(L)=\bigoplus_{w\mid v}H^1(L_w,E)(p),
 \end{equation}

 where $w$ runs over all primes of $L$ dividing $v$. We have the localisation map

 \begin{equation}
  \lambda_S(L):H^1(G_S(L),\Ep) \rightarrow \bigoplus_{v\in S}J_v(L).
 \end{equation}

 For an infinite algebraic extension $H$ of $F$, define
 \begin{equation}
  J_v(H)=\varinjlim J_v(L),
 \end{equation}

 where $L$ runs over all finite extensions of $F$ contained in $L$ and the inductive limit is taken with respect to the restriction maps. Also, define
 $\lambda_S(H)$ to be the inductive limit of the localisation maps $\lambda_S(L)$. It is well-known that the Selmer group Sel$(E/H)$ is given by the the
 exact sequence
 \begin{equation}
  0 \rightarrow \Sel(E/H) \rightarrow H^1(G_S(H),\Ep) \rightarrow \bigoplus_{v\in S}J_v(H).
 \end{equation}
 In particular, the Selmer group of $E$ over $F_\infty$ is given by
 \begin{equation}\label{Selmer*}
  0 \rightarrow \Sel(E/F_\infty) \rightarrow H^1(G_S(F_\infty),\Ep) \rTo^{\lambda_S(F_\infty)} \bigoplus_{v\in S}J_v(F_\infty).
 \end{equation}
 Taking $\Sigma$-invariants of ~\eqref{Selmer*}, we obtain the fundamental diagram
 \begin{diagram}\label{fundamentaldiag}
   0   &    \rTo & \SE^{\Sigma}                & \rTo & H^1(G_S(F_\infty),\Ep)^{\Sigma}           & \rTo^{\psi_S(F_\infty)}  &  \bigoplus_{v\in S}J_v(F_\infty)^{\Sigma}\\
       &         &   \uTo_{\alpha}             &      &   \uTo_{\beta}                            &       &    \uTo_{\gamma}                         \\
   0   &    \rTo & \Sel(E/F)                   & \rTo & H^1(G_S(F),\Ep)                           &\rTo   &  \bigoplus_{v\in S}J_v(F)                \\
 \end{diagram}
 This diagram is the basic tool for studying the Selmer group $\SE$ as a $\Lambda(\Sigma)$-module. Note that $\gamma$ is the direct sum
 \begin{equation}
  \gamma=\bigoplus_{v\in S}\gamma_v
 \end{equation}
 of the local restriction maps
 \begin{equation}
  \gamma_v:J_v(F) \rightarrow J_v(F_\infty)^ \Sigma.
 \end{equation}
 In Section 2 we study the kernels and cokernels of the local restriction maps - most of the results are immediate generalisations of the corresponding
 results in~\cite{coates2}. While $\lambda_S(L)$ is, in general, not surjective for finite extensions $L$ of $F$, it seems to be true (although we can only prove
 it in very few cases) that $\lambda_S(F_\infty)$ should always be surjective. Assuming the surjectivity of $\lambda_S(F_\infty)$, we use the fundamental diagram to calculate the $\Sigma$-Euler
 characteristic of $\SE$ in the case when $\Ep$ is rational over $F_\infty$ in Section 3. We then finish the proof of Theorem~\ref{EC} in this case by 
 relating the surjectivity of $\lambda_S(F_\infty)$ to the structure of ${\cal C}(E/F_\infty)$ as a $\Lambda(\Sigma)$-module. In Section 4 we explain
 how to adapt the arguments to the case when $\Ep$ is not rational over $F_\infty$. Finally, in Section 5 we illustrate our general theory with some numerical
 examples. 

\section{Local Results}\label{local}

 \subsection{Decomposition of primes in $F_\infty$}\label{decomp}

  As before, $E$ is an elliptic curve defined over a number field $F$ and $A$ is an abelian variety defined over
  $F$. Throughout this chapter, $p$ will be a prime $\geq 5$, $F_\infty=F(A_{p^\infty})$ and $\Sigma=$Gal$(F_\infty/F)$. We will assume that $\Sigma$
  has no element of order $p$.
  For a prime $v$ of $F$ and a fixed prime $w$ of $F_\infty$ above $v$, we let
  $\Sigma_w=$Gal$(F_{\infty,w}/F_v)$, which can be identified with the decomposition group $D(v)$ in $\Sigma$ of any fixed prime $w$ of $F_\infty$ above $v$.

  \begin{lem}\label{locdim}
   Assume that $v$ does not divide $p$. Then $D(v)$ has dimension $1$ or $2$, according as $A$ does or does not have potentially good reduction at $v$. 
  \end{lem}
  \begin{proof}
   We can identify $D(v)$ with the local Galois group $\Sigma_w$. Suppose that $A$ has 
   potentially good reduction at $v$. Since the dimension of $\Sigma_w$
   does not change when $F_v$ is replaced by some finite extension of $F_v$ contained in $F_{\infty,w}$, we may assume
   that $A$ has good reduction at $v$. Then the extension $F_{\infty,w}$ over $F_v$ is unramified and therefore generated by the Frobenius automorphism, so dim$(\Sigma_w)=1$. 
   Suppose that $A$ does not have potentially good reduction at $v$. 
   By the Weil pairing, $F_v^{cyc}\subset \Fiw$. Now dim$(F_v^{cyc}/F_v)=1$, so dim$(\Sigma_w)\geq 1$. 
   Since $v\nmid p$, the extension $F_v^{cyc}$ over $F_v$ is unramified. Denote by $I_v$ the inertia group of $\Sigma_w$. Local class field theory shows 
   that dim$(\Sigma_w)\leq 2$, so it is sufficient to prove that $I_v$ has positive $p$-cohomological dimension. By the criterion of 
   Neron-Ogg-Shafarevich, $I_v$ is infinite. Since $\Sigma_w$ contains a pro-$p$ subgroup of finite index, the same is true for $I_v$, so dim$(I_v)\geq 1$.
  \end{proof} 
  
 \subsection{Local cohomology}\label{LocCoh}

  \begin{lem}\label{Shapiro}
   Let $v \in S$. Then for all $i \geq 0$ there is a canonical isomorphism
   \begin{equation}
    H^i(\Sigma,J_v(F_\infty))\cong H^i(\Sigma_w,H^1(F_{\infty,w},E)(p)).
   \end{equation}
  \end{lem}
  \begin{proof}
   Easy consequence of Shapiro's lemma.
  \end{proof}


  \begin{prop}\label{locvanish}
   Let $v$ be a prime of $F$ not dividing $p$. Then
   \begin{equation}
    H^i(\Sigma,J_v(F_\infty))=0
   \end{equation}
   \noindent for all $i \geq 1$.
  \end{prop}

  \noindent The proposition will follow from Lemmas~\ref{locvanishbr} and~\ref{locvanishgr} below:-

  \begin{lem}\label{locvanishbr}
   If $v$ is a prime of $F$ dividing $p$ where $A$ does not have potentially good reduction, then
   \begin{equation}
    J_v(F_\infty)=0.
   \end{equation}   
  \end{lem}  
  \begin{proof}
   By local Kummer theory, $H^1(F_{\infty,w},E)(p) \cong H^1(F_{\infty,w},E_{p^\infty})$. Lemma~\ref{locdim} shows that
   $A$ does not have potentially good reduction reduction at $v$ if and only if dim$(\Sigma_w)=2$. But in this case the profinite degree of the extension
   $\bar{F}_v$ over $F_{\infty,w}$ is coprime to $p$, so cd$_p($Gal$(\bar{F}_v/F_{\infty,w}))=0$.
  \end{proof}   

  \begin{lem}\label{locvanishgr}
   If $v$ is a prime of $F$ not dividing $p$ where $A$ has potentially good reduction, then    
   \begin{equation}
    H^i(\Sigma,J_v(F_\infty))=0
   \end{equation}
   for all $i \geq 1$.
  \end{lem}
  \begin{proof}
   \noindent By Lemma ~\ref{Shapiro}, it is required to show that
   \begin{equation}\label{lv}
    H^i(\Sigma_w,J_v(F_\infty))=0
   \end{equation}
   \noindent for all $i \geq 1$. Since dim$(\Sigma_w)=1$ and $\Sigma_w$ has no element of order $p$, we have cd$_p(\Sigma_w)=1$, so ~\eqref{lv} is true
   for all $ i \geq 2$. 
   \noindent Note that by Tate local duality,
   \begin{equation}
    H^2(L,E_{p^\infty})=0
   \end{equation}
   \noindent for every finite extension $L$ of $F_v$. Allowing $L$ to range over all finite extensions of $F_v$ contained in $F_{\infty,w}$, we
   deduce that
   \begin{equation}
    H^2(F_{\infty,w},E_{p^\infty})=0.
   \end{equation}
   \noindent We conclude from the Hochschild-Serre spectral sequence that the sequence
   \begin{equation}
    H^2(F_v,E_{p^\infty}) \rightarrow H^1(\Sigma_w,H^1(F_{\infty,w},E_{p^\infty})) \rightarrow H^3(\Sigma_w,\Ep(\Fiw))
   \end{equation}
   \noindent is exact. But the groups on the left and right are zero which proves~\eqref{lv} for $i=1$.
  \end{proof}


  In the following, let $k_{\infty,w}$ denote the residue field of $\Fiw$.

  \begin{lem}\label{locisom}
   Let $v$ be a prime of $F$ dividing $p$ and assume that both $A$ and $E$ have good ordinary reduction at $v$. Then 
   \begin{equation}
    H^i(\Sigma_w,H^1(F_{\infty,w},E)(p)) \cong H^{i+2}(\Sigma_w,\tE(\kiw))
   \end{equation}
   \noindent for all $i \geq 1$.
  \end{lem}
  \begin{proof}
   Recall that $F_{\infty,w}$ is deeply ramified in the sense of ~\cite{coatesgreenberg} since it contains the deeply ramified field $F_v(\mu_{p^\infty})$,
   so there is a canonical $\Sigma_w$-isomorphism
   \begin{equation}\label{reducisom}
    H^1(F_{\infty,w},E)(p) \cong H^1(F_{\infty,w},\tE).
   \end{equation}
   Again using Tate local duality, we see that
   \begin{equation}
    H^i(F_{\infty,w},\tilde{E}_{v,p^\infty})=0
   \end{equation}
   \noindent for all $i \geq 2$. The lemma now follows from applying Hochschild-Serre for the module on the right of~\eqref{reducisom} to the extension 
   $F_{\infty,w}$ over $F_v$. 
  \end{proof}


 \subsection{Analysis of the local restriction maps}\label{analoc}
 
  Recall that for each prime $v$ of $S$ we have defined the local restriction map
  \begin{equation}
   \gamma_v:J_v(F)=H^1(F_v,E)(p) \rightarrow J_v(F_\infty)^{\Sigma}=H^1(F_{\infty,w},E)(p)^{\Sigma_w},
  \end{equation}
  \noindent where 
  \begin{equation}
   J_v(F_\infty)=\varinjlim\bigoplus_{w\mid v}H^1(F_{n,w},E)(p).
  \end{equation}

  We quote the following result from \cite{coates2}:-

  \begin{prop}\label{gammanondivp}
   Let $v$ be any finite prime of $F$ not dividing $p$. If both $A$ and $E$ have potentially good reduction at $v$, then $\gamma_v$ is an
   isomorphism. If $A$ has potentially good reduction at $v$ and $E$ does not have potentially good reduction at $v$, then $\gamma_v$ is
   surjective, and $\#\ker(\gamma_v)=\mid c_v\mid_p^{-1}$. If $A$ does not have potentially good reduction at $v$, then $\gamma_v$ is the zero map,
   and the order of its kernel is $\mid \frac{c_v}{L_v(E,1)}\mid _p^{-1}$.
  \end{prop}
  \begin{proof}
   We only give a sketch of the proof - for the details see Proposition 3.9 in~\cite{coates2}. 
   By inflation-restriction, we have
   \begin{align}
    \ker (\gamma_v)&=H^1(\Sigma_w,\Ep(\Fiw)),\\
    \coker (\gamma_v)&=H^2(\Sigma_w,\Ep(\Fiw)).
   \end{align}
   If $A$ does not have potentially good reduction at $v$, then $J_v(F_\infty)=0$ by Lemma~\ref{locvanishbr}, so
   \begin{equation}
    \ker(\gamma_v)=J_v(F_\infty),
   \end{equation}
   which is of order $\mid \frac{c_v}{L_v(E,1)}\mid _p^{-1}$. Assume now that $A$ has potentially good reduction at $v$. As shown in
   Lemma~\ref{decomp}, we have $\dim(\Sigma_w)=1$, so $\coker(\gamma_v)=0$. Put $\Gamma=\Gal(F_v^{cyc}/F_v)$ and $\Delta=\Gal(\Fiw/F_v^{cyc})$.
   Then we have the exact sequence
   \begin{equation}
    0\rTo H^1(\Gamma,\Ep(F_v^{cyc}))\rTo\ker(\gamma_v)\rTo H^1(\Delta,\Ep(\Fiw)).
   \end{equation}
   It is shown in~\cite{coates2} that the term on the right is zero and the term on the left is finite of order $\mid c_v\mid_p^{-1}$. If
   $\ord_v(j_E)\geq 0$, then as shown in~\cite{coates2} the only primes which divide $c_v$ lie in the set $\{ 2,3\}$, which proves the
   proposition.  
  \end{proof}

  For a prime $v$ of $F$ dividing $p$, the situation is more delicate; we again use the fact that the field $F_{\infty,w}$, for $w$ dividing $p$ is
  deeply ramified.

  \begin{prop}\label{gammadivp}
   Let $v$ be a prime of $F$ dividing $p$, and assume that both $A$ and $E$ have good ordinary reduction at $v$. Then both $\ker(\gamma_v)$ and 
   $\coker(\gamma_v)$ are finite, and
   \begin{equation}
    \frac{\#\ker(\gamma_v)}{\#\coker(\gamma_v)}=\frac{h_{v,0}h_{v,1}}{h_{v,2}},
   \end{equation}
  \noindent where
   \begin{equation}
    h_{v,i}=\#H^i(\Sigma_w,\tE(\kiw)).
   \end{equation}
  \end{prop}  

  The proof of Proposition~\ref{gammadivp} will be split up into a series of lemmas. The finiteness of the $H^i(\Sigma_w,\tE(\kiw))$ is shown in 
  Proposition~\ref{locEC} below. The proofs of Lemmas~\ref{m1}, \ref{m2}, \ref{m3} can be found in
  \cite{coates2}. In the following, let $\mathfrak{m}_\infty$ denote the maximal ideal in the ring of integers of $F_{\infty, w}$ and
  $\bar{\mathfrak{m}}$ the maximal ideal of the ring of integers of $\bar{F_v}$. Let $O_v$ be the ring of integers of $F_v$, and let $\hat{E}_v$ be the formal
  group defined over $O_v$ giving the kernel of reduction modulo $v$ on $E$. 

  \begin{lem}\label{m1}
   For all $i \geq 1$
   \begin{equation}
    H^i(\Sigma_w,\hat{E}_v(\mathfrak{m}_\infty)) \cong H^i(F_v,\hat{E}_v(\bar{\mathfrak{m}})).
   \end{equation}
  \end{lem}

  \begin{lem}\label{m2}
   Let $v$ be a prime of $F$ dividing $p$, suppose that $E$ has good ordinary reduction at $v$. Then
   $H^1(F_v,\hat{E}_v(\bar{\mathfrak{m}}))$ is finite of order $e_v$, and
   \begin{equation}
    H^i(F_v,\hat{E}_v(\bar{\mathfrak{m}}))=0
   \end{equation}
   \noindent for all $i \geq 2$.
  \end{lem}

  \begin{lem}\label{m3}
   Let $v$ be a prime of $F$ dividing $p$, suppose that both $A$ and $E$ have good ordinary reduction at $v$. Then
   \begin{equation}
    H^i(\Sigma_w,E(F_{\infty,w}))(p) \cong H^i(\Sigma_w,\tE(\kiw))
   \end{equation}
   \noindent for all $i \geq 2$, and we have an exact sequence
   \begin{equation}
    0 \rightarrow H^1(F_v,\hat{E}_v(\bar{\mathfrak{m}})) \rightarrow H^1(\Sigma_w,E(F_{\infty,w}))(p) \rightarrow
    H^1(\Sigma_w,\tE(\kiw)) \rightarrow 0.
   \end{equation}
  \end{lem}

  Applying the Hochschild-Serre spectral sequence to the extension $F_{\infty,w}$ over $F_v$ and recalling that
  $H^2(F_v,E_{p^\infty})=0$, we see that

  \begin{align}
   \ker(\gamma_v) &=H^1(\Sigma_w,E(F_{\infty,w}))(p) \\
   \coker(\gamma_v) &=H^2(\Sigma_w,E(F_{\infty,w}))(p).
  \end{align}

  Using the notation defined above, it follows from the previous lemmas that

  \begin{align}
   \# \ker(\gamma_v) &=h_0^vh_1^v \\
   \# \coker(\gamma_v) &=h_2^v,
  \end{align}

  \noindent which finishes the proof of Proposition~\ref{gammadivp}.
   The following result about the local Euler characteristic will be important:-

  \begin{prop}\label{locEC}
   Let $v$ be a prime of $F$ dividing $p$. Assume that both $A$ and $E$ have good ordinary reduction at $v$. Then
   $\tilde{E}_{v,p^\infty}(\kiw)$ has finite $\Sigma_w$-Euler characteristic, and
   \begin{equation}
    \chi(\Sigma_w,\tE(\kiw))=1.
   \end{equation}
  \end{prop}

  We split the proof of the proposition up into two cases:-

  \begin{lem}\label{locECcase1}
   Assume that $\Ep$ is rational over $F_{\infty,w}$. Under the assumptions of Proposition ~\ref{locEC}, $\tilde{E}_{v,p^\infty}$ has finite 
   $\Sigma_w$-Euler characteristic, and
   \begin{equation}
    \chi(\Sigma_w,\tilde{E}_{v,p^\infty})=1.
   \end{equation}
  \end{lem}
  \begin{proof}
   Let $G=$Gal$(F_{\infty,w}/F^{cyc}_w)$. By the usual argument (c.f. $\it{A.2.9}$ in~\cite{coatessujatha2}) it is enough to show that 
   $H^ i(G,\tilde{E}_{v,p^\infty})$ is finite for all $i \geq 0$. As $\Ep$ is rational over $F_\infty$, we can replace $A$
   by $A \times E$ without changing $F_\infty$. Hence we may assume that $E$ is an isogeny factor of $A$, so $V_{\tilde{E}}=T_p(\tilde{E_v}) \otimes \QQ_p$ is a subquotient of $V_A=T_p(A)
   \otimes \QQ_p$ in the sense of~\cite{coatessujathawintenberger}. The representation of $G$ on $V_A$ is faithful, so $G \subset $GL$(V_A)$. Write $L(G) \subset$End$(V_A)$ for the Lie algebra
   of $G$. Let $\overline{\QQ}_p$ be a fixed algebraic closure of $\QQ_p$. For any vector space $W$ over $\QQ_p$, write $W_{\overline{\QQ}_p}$ for $W \otimes
   _{\QQ_p} \overline{\QQ}_p$. 
   As shown in~\cite{coatessujathawintenberger}, $L(G)_{\overline{\QQ}_p} \subset$End$(V_{A,\overline{\QQ}_p})$ satisfies the strong Serre criterion. By
   Proposition 2.5 of~\cite{coatessujathawintenberger}, this implies that 
   \begin{equation}
    H^ i(G,V_{\tilde{E}})=0
   \end{equation}
   for all $i \geq 0$, and hence that $H^ i(G,\tilde{E}_{v,p^ \infty})$ is finite for all $i \geq 0$.
  \end{proof}

  \begin{lem}
   Assume that $\Ep$ is not rational over $F_{\infty,w}$. Again under the assumptions of Proposition ~\ref{locEC}, $\tE(\kiw)$ has finite 
   $\Sigma_w$-Euler characteristic, and
   \begin{equation}
    \chi(\Sigma_w,\tE(\kiw))=1.
   \end{equation}
  \end{lem}
  \begin{proof}
   Put $K_{\infty,w}=F_{\infty,w}(\Ep)$. Define the Galois groups $R_w=\Gal(K_{\infty,w}/F_v)$ and $\Omega=\Gal(K_{\infty,w}/F_{\infty,w})$. Note that $K_{\infty,w}$ is obtained by
   adjoining to $F_v$ the $p$-division points of the abelian variety $A\times E$. We first show that $R_w$ has no element of order $p$. Assume that $x\in R_w$ is
   of order $p$. Then $x$ cannot act trivially on both $\Ep$ and $A_{p^\infty}$. Hence the image of $x$ in $\Sigma_w$ and $\Gal(F_v(\Ep)/F_v)$ is
   nontrivial, so has order $p$. Note that $\Gal(F_v(\Ep)/F_v)$ has no element of order $p$ since $p\geq 5$. Since by assumption $\Sigma_w$ does not contain an
   element of order $p$, we get a contradiction. By Lemma~\ref{locECcase1} and the above remark, it therefore
   follows that $\tE$ has finite
   $R_w$-cohomology, and
   \begin{equation}
    \chi(R_w,\tE)=1.
   \end{equation}
   Since $\Ep$ is not rational over $F_{\infty,w}$, Proposition 3.15 in~\cite{coates2} shows that either $\Omega$ is a finite abelian group of order
   prime to $p$ or $\Omega\cong \ZZ_p\times\Delta$, where $\Delta$ is some finite
   abelian group of order prime to $p$. We will first deal with the easy case when $\Omega$ is finite abelian of order prime to $p$.
   Applying Hochschild-Serre to the extension $\Kiw$ over $\Fiw$ gives isomorphisms
   \begin{equation}
    H^i(\Sigma_w,\tE(\kiw))\cong H^i(R_w,\tE)
   \end{equation}
   for all $i\geq 0$, so
   \begin{equation}
    \chi(\Sigma_w,\tE(\kiw))=\chi(R_w,\tE)=1.
   \end{equation}
   Assume now that $\Omega\cong \ZZ_p\times\Delta$, where $\Delta$ is finite abelian of order prime to $p$. It is shown in~\cite{coates2} that the fixed
   field of $\Delta$ is a ramified extension of $\Fiw$. \\
   Assume first that $\Delta$ is trivial. Note that in this case $\tE$ is rational over $\Fiw$. Now cd$_p(\Omega)=1$, so applying Hochschild-Serre to
   the extension $K_{\infty,w}$ over $\Fiw$ gives the long exact sequence
   \begin{align*}
    0 &\rightarrow H^1(\Sigma_w,\tE) \rightarrow H^1(R_w,\tE) \rightarrow H^1(\Omega,\tE)^{\Sigma_w} \\
      &\rightarrow H^2(\Sigma_w,\tE) \rightarrow H^2(R_w,\tE) \rightarrow H^1(\Sigma_w,H^1(\Omega,\tE)) \rightarrow \dots
   \end{align*}
   In order to prove the lemma we must therefore show that $H^1(\Omega,\tE)$ has finite $\Sigma_w$-cohomology and that
   \begin{equation}\label{helpEC}
    \chi(\Sigma_w,H^1(\Omega,\tE))=1.
   \end{equation}
   Let $H_\infty$ be the maximal unramified extension of $F_v$ contained in $\Fiw$ and $G=\Gal(\Fiw/H_\infty)$. Now $\Gal(H_\infty/F_v)$ is isomorphic to the direct product of
   $\ZZ_p$ with some finite abelian group of order prime to $p$, so be the usual argument (c.f. $\it{A.2.9}$ in~\cite{coatessujatha2}) it is enough to show
   that $H^i(G,H^1(\Omega,\tE))$ is finite for all $i\geq 0$. \\
   Let $M_\infty=H_\infty(\mu_{p^\infty})$ and $\Gamma=\Gal(M_\infty/H_\infty)$. By the Weil pairing, $M_\infty\subset\Fiw$. 
   It is shown in~\cite{greenberg2} that (i) $\Gal(K_{\infty,w}/M_\infty)\cong \ZZ_p^{n+1}$ for some $n$ and (ii)
   $\Gamma$ acts on $\Gal(K_{\infty,w}/M_\infty)$ as the cyclotomic character $\chi$. To lighten notation, let $X=\Gal(\Fiw/M_\infty)$. 
   Choose a $\ZZ_p$-filtration of $X$, so
   \begin{equation}
    \{1\}=X_n\supset\dots\supset X_1\supset X_0=X,
   \end{equation}
   where $X_{i-1}/X_{i}\cong\ZZ_p$ for all $0<i\leq n$. Note that $X_i$ is a normal subgroup of $G$ for all $0\leq i\leq n$ thanks to property (ii)
   above. For $0<i\leq n$, let $G_i=G/X_i$. Observe that cd$_p(G_i)=i+1$ for all $i$.   
   
  \vspace{1ex}      
   {\it Claim} Let $1\leq j\leq n$. Let $A$ be a discrete $G_j$-module satisfying (i) $X/X_j$ acts trivially on $A$ and
   (ii) $A\cong(\QQ_p/\ZZ_p)(\chi^{-k})$ as a $\Gamma$-module for some $k\geq 1$.
   Then $H^i(G_j,A)$ is finite for all $i\geq 0$.
   \vspace{1ex}
   
   To prove the claim we are going to use repeatedly the following elementary observation: Let $\Gamma$ be a profinite group which is isomorphic to the direct product of $\ZZ_p$ with some
   finite abelian group of order prime to $p$. Let $A$ be a discrete $\Gamma$-module which isomorphic to $\QQ_p/\ZZ_p$ and on which $\Gamma$ acts via a nontrivial
   continuous homomorphism $\gamma:\Gamma\rightarrow \ZZ_p^\times$ whose image in $\ZZ_p^\times$ is infinite. Then $H^0(\Gamma,A)$ is finite and 
   $H^1(\Gamma,A)=0$.
   \vspace{1ex}
   
   {\it Proof of claim.} We procede by induction on $j$. \\
   $j=1$: Since cd$_p(X_0/X_1)=1$, Hochschild-Serre gives the exact sequence
   \begin{equation}
    0 \rightarrow H^1(\Gamma,A) \rightarrow H^1(G_1,A) \rightarrow H^1(X_0/X_1,A)^\Gamma \rightarrow 0
   \end{equation}
   and isomorphisms
   \begin{equation}
    H^i(G_1,A)\cong H^{i-1}(\Gamma,H^1(X_0/X_1,A))
   \end{equation}
   Since $\Gamma$ acts on $M_1$ via the cyclotomic character, we have an isomorphism of $\Gamma$-modules $H^1(X_0/X_1,A)\cong(\QQ_p/\ZZ_p)(\chi^{-(k+1)})$.
   Using the above
   observation and the fact that cd$_p(G_1)=2$ shows that $H^0(G_1,A)$ and $H^1(G_1,A)$ are finite and $H^i(G_1,A)=0$ for all $i\geq 2$.\\
   Let $1<j\leq n$, suppose that the claim holds for $j-1$. Now cd$_p(X_{j-1}/X_j)=1$, so again Hochschild-Serre gives an exact sequence
   \begin{align*}
    0 &\rightarrow H^1(G_{j-1},A) \rightarrow H^1(G_j,A) \rightarrow H^1(X_{j-1}/X_j,A)^{G_{j-1}}\\ 
      &\rightarrow\dots \\
      &\rightarrow H^j(G_{j-1},A) \rightarrow H^j(G_j,A) \rightarrow H^{j-1}(G_{j-1},H^1(X_{j-1}/X_j,A))\rightarrow 0
   \end{align*}
   By hypothesis, $H^i(G_{j-1},A)$ is finite for all $i\geq 0$. Now $X/X_{j-1}$ acts trivially on $H^1(X_{j-1}/X_j,A)$ and
   $H^1(X_{j-1}/X_j,A)\cong(\QQ_p/\ZZ_p)(\chi^{-(k+1)})$ as a $\Gamma$-module, so $H^1(X_{j-1}/X_j,A)$ satisfies conditions (i) and (ii) above. It follows again from the
   induction hypothesis that $H^i(G_{j-1},H^1(X_{j-1}/X_j,A))$ is finite for all $i\geq 0$. 
   \vspace{1ex}
   
   The finiteness of the $H^i(G,H^1(\Omega,\tE))$ now follows by applying this result for $j=n$ to the $G$-module $H^1(\Omega,\tE)$. \\
   \vspace{1ex}
   
   Assume now that $\Delta$ is not trivial. Let $\Liw$ be the fixed field of the $\ZZ_p$-component of $H$. It is shown in~\cite{coates2} that  
   $\Liw$ is an unramified extension of $\Fiw$. Let $H_\infty$ denote
   the maximal unramified extension of $F_v$ contained in $\Liw$, $S_w=\Gal(\Liw/F_v)$ and $G_w=\Gal(\Liw/H_\infty)$. By the same argument as above, $H^i(G_w,\tE)$ is finite
   for all $i\geq 0$, so $\tE$ has finite $S_w$-cohomology and 
   \begin{equation}
    \chi(S_w,\tE)=1.
   \end{equation}
   Now $\Delta$ is finite of order prime to $p$, so it has $p$-cohomological dimension equal to 0. Applying Hochschild-Serre to the extension $\Liw$ over
   $\Fiw$ therefore gives isomorphisms
   \begin{equation}
    H^i(\Sigma_w,\tE(\kiw))\cong H^i(S_w,\tE)
   \end{equation}
   for all $i\geq 0$, which shows that $\tE(\kiw)$ has finite $\Sigma_w$-cohomology and
   \begin{equation}
    \chi(\Sigma_w,\tE(\kiw))=1.
   \end{equation}  
  \end{proof}


 \section{Global Calculations when $\Ep(F_\infty)=\Ep$}\label{globalrational}
  In this section, again $E$ is an elliptic curve defined over a number field $F$, $A$ is an abelian variety defined over $F$, $p$ a prime $\geq 5$,
  $F_\infty=F(A_{p^\infty})$ and $\Sigma=\Gal(F_\infty/F)$. We shall assume thoughout that $\Sigma$ has no element of order $p$, so that its $p$-cohomological 
  dimension is equal to its dimension as a $p$-adic Lie group and all the results from section~\ref{local} are valid. Also, we assume that
  $E_{p^\infty}$ is rational over $F_\infty$. We use without comment the notation of Section 2.


  \subsection{Global cohomology}

   For the calculation of the $\Sigma$-Euler characteristic of $\Sel(E/F_\infty)$, we need the following results:-

   \begin{prop}\label{globEC}
    The $\Sigma$-module $E_{p^\infty}$ has finite $\Sigma$-Euler characteristic, and
    \begin{equation}
     \chi(\Sigma,E_{p^\infty})=1.
    \end{equation}
   \end{prop}
   \begin{proof}
    Let $B=E \times A$. Now $\Ep$ is rational over $F_\infty$, so $F_\infty=F(B_{p^\infty})$, and we have an exact sequence of $\Sigma$-modules 
    \begin{equation}\label{shortexseqabvar}
     0 \rightarrow \Ep \rightarrow B_{p^\infty} \rightarrow A_{p^\infty} \rightarrow 0.
    \end{equation}
    As shown in \cite{coatessujatha2}, both $A_{p^\infty}$ and $B_{p^\infty}$ have finite $\Sigma$-cohomology, and 
    \begin{equation}
     \chi(\Sigma,B_{p^\infty})=\chi(\Sigma,A_{p^\infty})=1.
    \end{equation}
    Taking $\Sigma$-cohomology of ~\eqref{shortexseqabvar} shows that $H^i(\Sigma,\Ep)$ is finite for all $i \geq 0$ and
    \begin{equation}
     \chi(\Sigma,\Ep)=1.
    \end{equation}
   \end{proof}
   
   Throughout the section, we let
   \begin{equation}
    h_i=\# H^i(\Sigma,\Ep).
   \end{equation}

   \begin{prop}\label{globvanish}
    \begin{equation}\label{globvanisheq}
     H^i(G_S(F_\infty),\Ep)=0
    \end{equation}
    for all $i \geq 2$.
   \end{prop}
   \begin{proof}
    cd$_p(G_S(F_\infty))=2$, so ~\eqref{globvanisheq} is true for for all $ i \geq 3$. Since $\Ep$ is rational over $F_\infty$, the case $i=2$ follows from an
    argument of Ochi as given in the proof of Theorem 2.10 in~\cite{coates2}.
   \end{proof}


  \subsection{Calculation of the Euler characteristic}\label{ECcalc}

   To avoid frequent repetitions, we introduce two hypotheses 
   \vspace{1ex}
   
   \noindent (R) both $A$ and $E$ have good ordinary reduction at all primes $v$ of $F$ dividing $p$\\
   \noindent (S) The map $\lambda_S(F_\infty):H^1(G_S(F_\infty),\Ep) \rightarrow \bigoplus_{v \in S} J_v(F_\infty)$ appearing 
   in \eqref{Selmer*} is surjective. \\
   \\
   We shall discuss in Subsection~\ref{Surj} various equivalent formulations of (S).
   We assume thoughout this section that $\Sel(E/F)$ is finite, i.e. that both $E(F)$ and $\Sha(E/F)(p)$ are finite.

   \begin{lem}\label{globisom}
    We have canonical isomorphisms
    \begin{equation}
     H^i(\Sigma, H^1(G_S(F_\infty),\Ep)) \cong H^{i+2}(\Sigma,\Ep)
    \end{equation}
    for all $i \geq 1$.
   \end{lem}
   \begin{proof}
    It is well-known that under the condition that $\Sel(E/F)$ is finite, we have
    \begin{equation}
     H^2(G_S(F),\Ep)=0.
    \end{equation}
    The Lemma follows from Proposition~\ref{globvanish} and from applying Hochschild-Serre to the extension $F_\infty$ over $F$.
   \end{proof}

   We now have all the necessary information to prove the $\Sigma$-Euler characteristic formula for $\Sel(E/F_\infty)$ in a special case:

   \begin{prop}\label{specialEC}
    Assume that $\Sel(E/F)$ is finite and that (R) and (S) hold. Then
    $\SE$ has finite $\Sigma$-Euler characteristic, and 
    \begin{equation}
     \chi(\Sigma,\SE)=\rho_p(E/F) \times\mid\prod_{v \in \frak{M}} L_v(E,1)\mid _p.
    \end{equation}
   \end{prop}

   The proof of the proposition will be split up into a series of lemmas. 


   Let $d$ denote the dimension of $\Sigma$ as a $p$-adic Lie group. Recall that cd$_p(\Sigma)=d$ since $\Sigma$ has no elements of order $p$.
   Let $T$ be the set of those primes in $S$, not dividing $p$, where $A$ has potentially good reduction.

   \begin{lem}\label{cokerfinite}
    Suppose that assumption (R) holds. Then $\coker(\psi_S(F_\infty))$
    is finite, where $\psi_S(F_\infty)$ is the map in the fundamental diagram.
   \end{lem}
   \begin{proof}
    Consider the following commutative diagram with exact rows which is extracted from the fundamental diagram:- 
    \begin{diagram}\label{extract26666}
      0   &    \rTo & \text{Im}(\psi_S(F_\infty)) & \rTo & \bigoplus_{v \in S}J_v(F_\infty)^{\Sigma} & \rTo &         \coker(\psi_S(F_\infty)) & \rTo & 0 \\
          &         &   \uTo_{\delta}             &      &  \uTo_{\gamma}                            &       &    \uTo_{\epsilon}              &      &\\
      0   &    \rTo & \text{Im}(\lambda_S(F))     & \rTo & \bigoplus_{v\in S}J_v(F)                  &\rTo  &          \coker(\lambda_S(F))    & \rTo & 0\\
    \end{diagram}
    Here $\delta$ and $\epsilon$ are
    the obvious induced maps. It is well-known that the finiteness of $\Sel(E/F)$ implies that $\coker(\lambda_S(F))$ is finite of order $\Ep(F)$. Also, we have already seen
    in Section~\ref{local} that $\gamma$ has finite kernel and cokernel, so the finiteness of $\coker(\psi_S(F_\infty))$ follows by applying the snake lemma
    to above diagram. 
   \end{proof}
   
   \begin{lem}\label{Sel0}
    Under the conditions of Proposition~\ref{specialEC}, $H^0(\Sigma,\SE)$ is finite, and
    \begin{equation}
     \#H^0(\Sigma,\SE)=\xi_p(E/F)\#\coker(\psi_S(F_\infty))\frac{h_0h_2}{h_1}\prod_{v\mid p}\frac{h_{v,1}}{h_{v,0}h_{v,2}},
    \end{equation}
    where
    \begin{equation}
     \xi_p(E/F)=\rho_p(E/F) \times\mid\prod_{v \in \frak{M}} L_v(E,1)\mid _p.
    \end{equation}
   \end{lem}
   \begin{proof}
    Using the fact that $H^2(G_S(F),\Ep)=0$, the inflation-restriction exact sequence shows that
    \begin{align*}
     \ker(\beta) &=H^1(\Sigma,\Ep), \\
     \coker(\beta) &= H^2(\Sigma,\Ep),
    \end{align*}
    so by Proposition~\ref{globEC}
    \begin{equation}\label{beta}
     \frac{\#\ker(\beta)}{\#\coker(\beta)}=\frac{h_1}{h_2}.
    \end{equation}
    Combining Propositions~\ref{gammanondivp} and \ref{gammadivp} shows that
    \begin{equation}\label{gamma}
     \frac{\#\ker(\gamma)}{\#\coker(\gamma)}=\prod_{v\mid p}\frac{h_{v,0}h_{v,1}}{h_{v,2}}\times \mid\prod_{v \in \frak{M}}\frac{c_v}{L_v(E,1)}\mid _p^{-1},
    \end{equation}
    where $e_v=\#\tE(k_v)$. Applying the snake lemma to the commutative diagram in Lemma 12 shows that $\delta$ has finite kernel and cokernel, and
    \begin{equation}\label{delta}
     \frac{\#\ker(\delta)}{\#\coker(\delta)}=\frac{\#\ker(\gamma)}{\#\coker(\gamma)}\frac{\coker(\psi_S(F_\infty))}{h_0}.
    \end{equation}
    We also have the commutative diagram with exact rows:-
    \begin{diagram}\label{extract2}
      0   &    \rTo & \SE^{\Sigma}                & \rTo & H^1(G_S(F_\infty),\Ep)^{\Sigma}           & \rTo  & \text{Im}\psi_S(F_\infty)        & \rTo & 0 \\
          &         &   \uTo_{\alpha}             &      &  \uTo_{\beta}                             &       &    \uTo_{\delta  }               &      &\\
      0   &    \rTo & \text{Sel}(E/F)             & \rTo & H^1(G_S(F),\Ep)                           &\rTo   & \text{Im}\lambda_S(F_\infty)     & \rTo & 0\\
    \end{diagram}
    Applying the snake lemma to this diagram shows that $\SE^\Sigma$ is finite of order
    \begin{equation}\label{Selmer1}
     \#H^0(\Sigma,\SE)=\#\Sel(E/F)\frac{\#\coker(\beta)}{\#\ker(\beta)}\frac{\#\ker(\delta)}{\#\coker(\delta)}.
    \end{equation}
    Since $\Sel(E/F)$ is finite, we have $\Sel(E/F)=\Sha(p)$. Also, it is well-known that $c_v\leq 4$ if $v\notin \frak{M}$. Combining~\eqref{beta},
    \eqref{gamma}, \eqref{delta} and ~\eqref{Selmer1} proves the lemma.
   \end{proof}
  
   To lighten notation, we define
   \begin{equation}
    W_{v,i}=H^i(\Sigma_w,H^1(\Fiw,\tE))
   \end{equation}
   for $i\geq 0$, and we put
   \begin{equation}
    B_\infty=H^1(G_S(F_\infty),\Ep).
   \end{equation}
  
   \begin{lem}\label{Sel1}
    Under the conditions of Proposition~\ref{specialEC}, $\SE$ has finite $\Sigma$-Euler characteristic, and
    \begin{equation}
     \chi(\Sigma,\SE)=\frac{\#H^0(\Sigma,\SE)}{\#\coker(\psi_S(F_ \infty))}\frac{h_1}{h_0h_2}\prod_{v\mid p}\frac{h_{v,0}h_{v,2}}{h_{v,1}}.
    \end{equation}
   \end{lem}
   \begin{proof}
    Since (S) holds, we have the short exact sequence of $\Sigma$-modules
    \begin{equation}
     0 \rightarrow \SE \rightarrow H^1(G_S(F_\infty),\Ep) \rightarrow \bigoplus_{v \in S}J_v(F_\infty) \rightarrow 0.
    \end{equation}
   Taking $\Sigma$-cohomology and recalling Lemmas~\ref{locvanishbr} and~\ref{locvanishgr} gives rise to the long exact sequence
   \begin{align*}
    0 &\rightarrow \SE^\Sigma \rightarrow B_\infty^\Sigma \rightarrow \bigoplus_{v \in
    T}H^1(F_{\infty,w},\Ep)^{\Sigma_w}\bigoplus_{v\mid p}W_{v,0} \\
    &\rightarrow H^1(\Sigma,\SE) \rightarrow H^1(\Sigma,B_\infty) \rightarrow \bigoplus_{v\mid p}W_{v,1} \\
    &\rightarrow \dots \\
    &\rightarrow H^d(\Sigma,\SE) \rightarrow H^d(\Sigma,B_\infty) \rightarrow \bigoplus_{v\mid p}W_{v,d}
    \rightarrow 0.
   \end{align*}
   Using Lemma~\ref{globisom}, the sequence can be written as
   \begin{align*}
    0 &\rightarrow \SE^\Sigma \rightarrow B_\infty^\Sigma \rightarrow \bigoplus_{v \in
    T}H^1(F_{\infty,w},\Ep)^{\Sigma_w}\bigoplus_{v\mid p}W_{v,0} \\
    &\rightarrow H^1(\Sigma,\SE) \rightarrow H^3(\Sigma,\Ep) \rightarrow \bigoplus_{v\mid p}H^3(\Sigma_w,\tE) \\
    &\rightarrow \dots \\
    &\rightarrow H^d(\Sigma,\SE) \rightarrow H^{d+2}(\Sigma,\Ep) \rightarrow \bigoplus_{v\mid p}H^{d+2}(\Sigma_w,\tE)
    \rightarrow 0.
   \end{align*}
   Recall that cd$_p(\Sigma)=d$, cd$_p(\Sigma_w) \leq d$ for all primes $w$ of $F_\infty$ dividing $p$. We have shown in Propositions~\ref{locEC} and~\ref{globEC} that
   $H^i(\Sigma,\Ep)$ and $H^i(\Sigma_w,\tE)$ are finite for all $i$. It follows that $H^d(\Sigma,\Ep)$ and $H^d(\Sigma_w,\tE)$ are both finite and
   $p$-divisible, so they must be zero. It is therefore immediate from above exact sequence that 
   \begin{equation}
    H^{d-1}(\Sigma,\SE)=H^d(\Sigma,\SE)=0,
   \end{equation}
   and we can extract the exact sequence
   \begin{align*}
    0 &\rightarrow \coker(\psi_S(F_\infty)) \rightarrow H^1(\Sigma,\SE) \rightarrow H^3(\Sigma,\Ep) \rightarrow \dots \\
    &\rightarrow H^{d-2}(\Sigma,\SE) \rightarrow 0.
   \end{align*}
   It now follows that all the terms in this exact sequence are finite. As shown in Lemma~\ref{Sel0}, $\SE^\Sigma$ is finite. Comparing the orders
   of the groups in above exact sequence and using Propositions~\ref{locEC} and~\ref{globEC} proves the lemma.
  \end{proof}

  Combining Lemmas~\ref{Sel0} and \ref{Sel1} gives the formula for $\chi(\Sigma,\SE)$. This finishes the proof of Proposition~\ref{specialEC}.


 \subsection{Surjectivity of $\lambda_S(F_\infty)$}\label{Surj}

  For $m\geq 0$, let $F_m=F(A_{p^{m+1}})$. Recall that $R=$Gal$(F_\infty/F_0)$. In this section, we are going to relate the surjectivity of 
  $\lambda_S(F_\infty)$ to the structure of ${\cal C}(E/F_\infty)$
  as a $\Lambda(R)$-module. We will assume throughout that $\Sel(E/F)$ is finite and that $A$ has good ordinary reduction at all primes $v$ of $F$ dividing $p$.\\
  \vspace{1ex}
  For each prime $v$ of $F$ dividing $p$, define
   \begin{equation}
    \tau_v(E/F)=
     \begin{cases}
      [F_v:\QQ_p]& \text{if $E$ has potentially supersingular reduction at $v$} \\
      0& \text{otherwise}
     \end{cases}
   \end{equation}
  \noindent and let
  \begin{equation}
   \tau_p(E/F)=\sum_{v\mid p}\tau_v(E/F).
  \end{equation}

  \noindent The main result of this section is the following:

  \begin{thm}\label{surj}
   Let $p$ be a prime $\geq 5$, suppose that $\Sigma$ has no element of order $p$. Then the $\Lambda(R)$-rank of ${\cal C}(E/F_\infty)$ equals $[\Sigma:R]\tau_p(E/F)$ if and only if 
   $\lambda_S(F_\infty)$ is
   surjective.
  \end{thm}

  The idea of the proof is to compute the $\Lambda(R)$-coranks of the terms of Cassels' variant of the Poitou-Tate exact sequence
  \begin{equation}\label{PoitouTate}
   0 \rightarrow \Sel(E/F_\infty) \rightarrow H^1(G_S(F_\infty),\Ep) \rightarrow \bigoplus_{v \in S}J_v(F_\infty) \rightarrow
   \widehat{R(E/F_\infty)} \rightarrow 0.
  \end{equation}
  Here, $R(E/F_\infty)$ denotes the compact Selmer group which is defined as follows:- Recall that $\Sel(E/F_m,p^n)$ is defined by the
  exactness of
  \begin{equation}
   0 \rightarrow \Sel(E/F_m,p^n) \rightarrow H^1(G_S(F_m),E_{p^n}) \rightarrow \bigoplus_{w_m}H^1(F_{m,w_m},E),
  \end{equation}
  where the $w_m$ range over all non-archimedean places of $F_m$. Then the compact Selmer group over $F_m$ is defined as the projective limit
  \begin{equation}
   R(E/F_m)=\varprojlim_n \Sel(E/F_m,p^n).
  \end{equation}
  Passing to the projective limit over $m$ with respect to the corestriction maps gives the compact Selmer group of $E$ over $F_\infty$:-
  \begin{equation}
   R(E/F_\infty)=\varprojlim_m R(E/F_m).
  \end{equation}
  Note that the map~$\bigoplus_{v \in S}J_v(F_\infty) \rightarrow\widehat{R(E/F_\infty)}$ in~\eqref{PoitouTate} is surjective since 
  \begin{equation}
   H^2(G_S(F_\infty),\Ep)=0
  \end{equation}
  by Proposition~\ref{globvanish}.
  \vspace{1ex}
  In the next two propositions we calculate the $\Lambda(R)$-coranks of the two middle terms of~\eqref{PoitouTate}.  

  \begin{prop}\label{globcor}
   \begin{equation}
    \crk_{\Lambda(R)}H^1(G_S(F_\infty),\Ep)=[\Sigma:R][F:\QQ]
   \end{equation}
  \end{prop}
  \begin{proof}
  Here $R$ is a pro-$p$ $p$-adic Lie group with no element of order $p$, so by Theorem (1.1) of \cite{howson1} 
  \begin{equation}
   \crk_{\Lambda(R)}H^1(G_S(F_\infty),\Ep)=\sum_{i \geq 0}(-1)^i\rk_{\ZZ_p}H^i(R,H^1(G_S(F_\infty),\Ep)). \notag
  \end{equation}
  Let $K=F_\infty^R$. By Theorem~\ref{globvanish} we have $H^i(G_S(F_\infty),\Ep)=0$ for all $i \geq 2$, so we can apply Hochschild-Serre to the extension $F_\infty$
  over $K$ to get the exact sequence
  \begin{gather*}
   0 \rightarrow H^1(R,\Ep) \rightarrow H^1(G_S(K),\Ep) \rightarrow H^1(G_S(F_\infty),\Ep)^\Sigma \rightarrow H^2(R,\Ep)\\
    \rightarrow H^2(G_S(K),\Ep) \rightarrow H^1(R,H^1(G_S(F_\infty),\Ep)) \rightarrow H^3(R,\Ep) \rightarrow 0
  \end{gather*}
  and isomorphisms
  \begin{equation}
   H^i(R,H^1(G_S(F_\infty),\Ep)) \cong H^{i+2}(R,\Ep)
  \end{equation}
  for all $i \geq 2$. By ~\eqref{globEC}, $H^i(R,\Ep)$ is finite for all $i \geq 0$, so
  \begin{equation*}
   \rk_{\ZZ_p}H^i(R,H^1(G_S(F_\infty),\Ep))=
    \begin{cases}
     \rk_{\ZZ_p}H^{i+1}(G_S(K),\Ep)& \text{for $i=0,1$}\\
     0& \text{for $i \geq 2$}                                        
    \end{cases}
  \end{equation*}           
  Let $K_\infty=K(\Ep)$. Then $G=$Gal$(K_\infty/K)$ is a pro-$p$ $p$-adic Lie group with no element of order $p$ since $p\geq 5$. Repeating above argument for the
  extension $K_\infty$ over $K$ shows that
  \begin{equation}
   \crk_{\Lambda(G)}H^1(G_S(K_\infty),\Ep)=\rk_{\ZZ_p}H^2(G_S(K),\Ep)-\rk_{\ZZ_p}H^1(G_S(K),\Ep), \notag
  \end{equation}
  so
  \begin{equation}
   \crk_{\Lambda(R)}H^1(G_S(F_\infty),\Ep)=\crk_{\Lambda(G)}H^1(G_S(K_\infty),\Ep). \notag
  \end{equation}
  In \cite{howson3}, Howson has shown that 
  \begin{equation*}
   \crk_{\Lambda(G)}H^1(G_S(K_\infty),\Ep)=[K:\QQ],
  \end{equation*} 
  which is equal to $[K:F][F:\QQ]=[\Sigma:R][F:\QQ]$.
 \end{proof}

 \begin{prop}\label{loccor}
  If $p \geq 5$, then 
  \begin{equation}
   \crk_{\Lambda(R)}\bigoplus_{v \in S}J_v(F_\infty)=[\Sigma:R]([F:\QQ]-\tau_p(E/F)).
  \end{equation}
 \end{prop}
 \begin{proof}
  As in the proof of \ref{globcor}, we use Theorem 1.1 of \cite{howson1}:
  \begin{align*}
   \crk_{\Lambda(R)}\bigoplus_{v \in S}J_v(F_\infty) &=\bigoplus_{v \in S}\crk_{\Lambda(R)}J_v(F_\infty)\\
   &=\sum_{v \in S}\sum_{i \geq 0}(-1)^i\rk_{\ZZ_p}H^i(R,J_v(F_\infty)).
  \end{align*}
  We analyse the primes $v$ in $S$ separately:
  \vspace{1ex}
  \begin{enumerate}
   \item If $v \nmid p$ and $A$ does not have potentially good reduction at $v$:-\\
    \vspace{1ex}
    \noindent We have $J_v(F_\infty)=0$ by \ref{locvanishbr}, so
    \begin{equation}
     \crk_{\Lambda(R)}J_v(F_\infty)=0.
    \end{equation}
   \item If $v \nmid p$ and $A$ has potentially good reduction at $v$:-\\
    \vspace{1ex}
     Lemma~\ref{locvanishgr} shows that $H^i(R,J_v(F_\infty))=0$ for all $i \geq 0$, so
    \begin{equation}
     \crk_{\Lambda(R)}J_v(F_\infty)=\rk_{\ZZ_p}J_v(F_\infty)^R.
    \end{equation}
    But as shown in ~\cite{howson3}, $J_v(F_\infty)^R$ is finite, so
    \begin{equation}
     \crk_{\Lambda(R)}J_v(F_\infty)=0.
    \end{equation}
   \vspace{1ex}
   
   Again, let $K=F_\infty^R$. Note that by Shapiro's lemma we have isomorphisms
   \begin{equation}
    H^i(R,J_v(F_\infty)) \cong \bigoplus_{v'\mid v,v' \text{ in } K}H^i(R_w,H^1(F_{\infty,w},E)(p))
   \end{equation}
   for any fixed choice of $w$ in $F_\infty$ above each $v'$.
   \item If $v \mid p$ and $E$ has good supersingular reduction at $v$:-\\
    \vspace{1ex}
    Since $F_{\infty,w}$ is a deeply ramified extension of $F_v$, Proposition (4.8) of \cite{coatesgreenberg} applies to show that
    $H^1(F_{\infty,w},E)(p)=0$, so
    \begin{equation}
     \crk_{\Lambda(R)}J_v(F_\infty)=0.
    \end{equation}
   \item If $v\mid p$ and $E$ has good ordinary reduction at $v$:-\\
   \vspace{1ex}
    We have
    \begin{equation}
     \crk_{\Lambda(R)}J_v(F_\infty)=\sum_{v' \mid v\text{ in }K}\sum_{i \geq
     0}\rk_{\ZZ_p}H^i(R_w,H^1(F_{\infty,w},\tEE)).
    \end{equation}
    Let $v'$ be a prime of $K$ dividing $v$. Applying Hochschild-Serre to the extension $\Fiw$ over $K_{v'}$ and using Tate local duality
    gives the exact sequence
    \begin{align*}
     0 &\rightarrow H^1(R_w,\tEE) \rightarrow H^1(K_{v'},\tEE) \rightarrow H^1(\Fiw,\tEE)^{R_w} \\
       &\rightarrow H^2(R_w,\tEE) \rightarrow 0
    \end{align*}
    and isomorphisms
    \begin{equation*}
     H^i(R_w,H^1(\Fiw,\tEE)) \cong H^{i+2}(R_w,\tEE)
    \end{equation*}
    for all $i \geq 2$. As shown in~\cite{howson1}, $H^i(R_w,\tEE)$ is finite for all $i \geq 0$, so 
    \begin{equation*}
     \rk_{\ZZ_p}H^i(R_w,H^1(\Fiw,\tEE))=
     \begin{cases}
      \rk_{\ZZ_p}H^1(K_{v'},\tEE)& \text{if $i=0$}\\
                                0& \text{if $i \geq 1$}
     \end{cases}
    \end{equation*}
    Let $K_\infty=K(E_{p^\infty})$. Again, $G=$Gal$(K_\infty/F)$ is a pro-$p$ $p$-adic Lie group with no elements of order $p$ since $p\geq 5$. Repeating above
    argument for the extension $K_\infty$ over $K$ shows that
    \begin{equation}
     \crk_{\Lambda(G)}J_{v'}(K_\infty)=\rk_{\ZZ_p}H^1(K_{v'},\tEE). \notag
    \end{equation}
    Howson has shown that 
    \begin{equation}
     \crk_{\Lambda(G)}J_{v'}(K_\infty)=[K_{v'}:\QQ_p], \notag
    \end{equation}
    so
    \begin{align*}
     \crk_{\Lambda(R)}J_v(F_\infty) &=\sum_{v'\mid v\text{ in }K}[K_{v'}:\QQ_p]\\
                                             &=[F_v:\QQ_p]\sum_{v'\mid v\text{ in }K}[K_{v'}:F]\\
					     & =[F_v:\QQ_p][K:F]\\
					     &=[F_v:\QQ_p][\Sigma:R].
    \end{align*}
  \end{enumerate}
  Combining results (1)-(4) shows that ideed
  \begin{equation}
   \crk_{\Lambda(R)}\bigoplus_{v \in S}J_v(F_\infty)=[\Sigma:R]([F:\QQ]-\tau_p(E/F)). \notag
  \end{equation}
 \end{proof}

 {\it Proof of \ref{surj}.} The implication that $\lambda_S(F_\infty)$ is surjective implies that ${\cal C}(E/F_\infty)$ has $\Lambda(R)$-rank
 $[\Sigma:R]\tau_p(E/F)$ follows immediately from Propositions~\ref{globcor} and \ref{loccor}. Conversely, suppose that ${\cal C}(E/F_\infty)$ has the expected
 $\Lambda(R)$-rank. The above propositions show that the dual of coker$(\lambda_S(F_\infty))$ is $\Lambda(R)$-torsion. The exact sequence ~\eqref{PoitouTate}
 shows that coker$(\lambda_S(F_\infty)) \cong \widehat{R(E/F_\infty)}$. However, Proposition (5.30) in \cite{howson3} generalizes immediately to show
 that $\widehat{R(E/F_\infty)}$ is $\Lambda(R)$-torsion free.

 \begin{cor}\label{surjtors}
  If $\Sigma$ has no element of order $p$ and ${\cal C}(E/F_\infty)$ is $\Lambda(R)$-torsion, then $\lambda_S(F_\infty)$ is surjective.
 \end{cor}

 \begin{rem} Propositions~\ref{globcor} and \ref{loccor} also give the following inequalities for the $\Lambda(R)$-rank of ${\cal C}(E/F_\infty)$,
 which is again analogous to the result in the GL$_2$-case discussed in \cite{coates2}:- 

 \begin{cor}\label{ineqcor}
  For all primes $p > 2$dim$(A)+1$,
  \begin{equation}
   \tau_p(E/F) \leq \frac{\rk_{\Lambda(R)} {\cal C}(E/F_\infty)}{[\Sigma:R]} \leq [F:\QQ].
  \end{equation}
 \end{cor}
 \end{rem}

 We conjecture that in fact

 \begin{con}
  \begin{equation}
   \rk_{\Lambda(R)}{\cal C}(E/F_\infty)=[\Sigma:R]\tau_p(E/F),
  \end{equation}
  which is a natural generalisation of a conjecture in~\cite{schneider}.
 \end{con}


 \section{Global Calculations when $\Ep(F_\infty)\neq \Ep$}\label{globalrational*}
  In this section, again $E$ is an elliptic curve defined over a number field $F$. I am grateful to Susan Howson for suggesting that the following result might
  be true:- 

  \begin{prop}
   Let $p$ be an odd prime and let $F_\infty$ be a $p$-adic Lie extension which contains the cyclotomic $\ZZ_p$-extension. Suppose that $\Ep$ is not rational 
   over $F_\infty$. Then $\Ep(F_\infty)$ is finite.
  \end{prop}
  \begin{proof}
   Suppose that $\Ep(F_\infty)$ is infinite. Assume first that $E$ does not have complex multiplication. If $\Ep(F_\infty)$ has
   $\ZZ_p$-corank equal to $1$, then $V_p(E)=T_p(E)\otimes_{\ZZ_p}\QQ_p$ has a $1$-dimensional $\Gal(\bar{F}/F)$-invariant subspace. However,  
   as shown in~\cite{serre5}, the Galois group of $F(\Ep)$ over $F$ is an open subgroup of GL$_2(\ZZ_p)$, so the representation of 
   $\Gal(\bar{F}/F)$ on $V_p(E)$ is irreducible, which gives the required contradiction.
   Assume now that $E$ has complex multiplication. We first show that the $p$-torsion points of $E$ are rational over $F_\infty$. Let 
   $\Delta=\Gal(F_\infty(E_p)/F_\infty)$. It is well-known that the order of $\Delta$ is prime to $p$ (cf.~\cite{perrin-riou}). By assumption, 
   $E_p(F_\infty)\neq 0$. 
   By choosing a suitable $\ZZ/p\ZZ$-basis of $E_p$, we get an isomorphism of $\Delta$ with a subgroup of GL$_2(\ZZ/p\ZZ)$ consisting of matrices of the form 
   $\begin{pmatrix}1 &x \\ 0 &1\end{pmatrix}$, where $x\in\ZZ/p\ZZ$. But if $x\neq 0$, then $\begin{pmatrix}1 &x \\ 0 &1\end{pmatrix}$ generates a 
   subgroup of order $p$. It follows that $\Delta$ is trivial and hence $E_p$ is rational over $F_\infty$. It is shown in~\cite{perrin-riou} that
   $\Gal(F(\Ep)/F)\cong \ZZ_p\times\ZZ_p$. Note that $F(\Ep)$ contains $F^{cyc}$ by the Weil pairing. Let $K_\infty$ denote the extension of $F$ generated 
   by $\Ep(F_\infty)$. Then $\Gal(K_\infty/F)\cong\ZZ_p\times\ZZ/p^n\ZZ$ for some $n\geq 0$. It is shown in~\cite{coatessujatha1} 
   that $\Ep(F^{cyc})$ is finite, so
   $K_\infty$ and $F^{cyc}$ intersect in a finite extension of $F$. It follows that $F(\Ep)=K_\infty F^{cyc}$, which proves the proposition.
  \end{proof}

  For the rest of this section, we let $p$ be a prime $\geq 5$, $A$ an abelian variety defined over $F$,
  $F_\infty=F(A_{p^\infty})$ and $\Sigma=\Gal(F_\infty/F)$. We shall assume thoughout that $\Sigma$ has no element of order $p$, so that its 
  $p$-cohomological dimension
  equals its dimension as a $p$-adic Lie group and all the results from Section~\ref{local} are valid. We now assume that
  $\Ep(F_\infty)\neq\Ep$.
  \vspace{1ex}

  We need the following result about the $\Sigma$-Euler characteristic of $\Ep(F_\infty)$:-

  \begin{prop}\label{globEC*}
   The $\Sigma$-module $\Ep(F_\infty)$ has finite $\Sigma$-Euler characteristic, and
   \begin{equation}
    \chi(\Sigma,\Ep(F_\infty))=1.
   \end{equation}
  \end{prop}
  \begin{proof}
   We are going to use the following result from~\cite{serre1}: Let $G$ be a pro-$p$ $p$-adic Lie group, let $M$ be a finite $p$-primary discrete $G$-module.
   Then $H^i(G,M)$ is finite for all $i$. 
   Since $\Ep(F_\infty)$ is finite, this proves the proposition when $\Sigma$ is pro-$p$. Suppose that $\Sigma$ is not pro-$p$. Let
   $G=\Gal(F_\infty/F^{cyc})$. $G$ is a closed subgroup of $\Sigma$, so it is itself a $p$-adic Lie group. As before, the proposition will follow if we can show 
   that $H^i(G,\Ep(F_\infty))$ is finite for all $i$. Fix an open normal pro-$p$ subgroup $H$ of
   $G$. Since $H$ is pro-$p$, $H^i(H,\Ep(F_\infty))$ is finite for all $i\geq 0$. The finiteness of the $H^i(G,\Ep(F_\infty))$ now follows from the 
   Hochschild-Serre spectral sequence
   \begin{equation}
    H^p(G,H^q(H,\Ep(F_\infty)))\Rightarrow H^{p+q}(G,\Ep(F_\infty)).
   \end{equation}
  \end{proof}
 
 Recall that $R=$Gal$(F_\infty/F_0)$. As an analogue of Proposition~\ref{globvanish}, we quote the following result from~\cite{hachimorivenjakob}:

 \begin{prop}
  If ${\cal C}(E/F_\infty)$ is $\Lambda(R)$-torsion, then \\
  (i) $H^2(G_S(F_\infty),\Ep)=0$ and \\
  (ii) The map
       \begin{equation}
        \lambda_S(F_\infty):H^1(G_S(F_\infty),\Ep)\rightarrow\bigoplus_{v\in S}J_v(F_\infty)
       \end{equation}
       is surjective.
 \end{prop}
 
 Assuming that ${\cal C}(E/F_\infty)$ is $\Lambda(R)$-torsion, the proof of Theorem~\ref{EC} now follows as described in Subsection~\ref{ECcalc}.


 \section{Numerical Examples}\label{Numerical}

 Under the circumstances of Theorem~\ref{EC}, it is in general difficult to check whether condition (iv) is satisfied. However, the following result,
 which is an immediate generalisation of Theorem 6.4 in~\cite{coateshowson2},
 turns out to be very useful in practice:

 \begin{thm}\label{finitegeneration}
  Let $A$ be an abelian variety defined over a finite extension $F$ of $\QQ$, let $p$ be a prime $\geq 5$, let
  $F_\infty = F(A_{p^{\infty}})$ and $\Sigma=\Gal(F_\infty/F)$. Let $E$ be an elliptic curve defined over $F$.
  Assume that there exists a finite extension $L$ of $F$ contained in $F_\infty$ such that $\Gal(F_\infty/L)$ is pro-$p$ and ${\cal
  C}(E/L^{cyc})$ is a finitely generated $\ZZ_p$-module. Then ${\cal C}(E/F_\infty)$ is a finitely generated $\Lambda(H)$-module. If $\Sigma$ has no element of
  order $p$, then it is $\Lambda(\Sigma)$-torsion.
 \end{thm}
 \begin{proof}
  It is easy to see that ${\cal C}(E/F_\infty)$ is finitely generated as a $\Lambda(\Sigma)$-module (c.f. Theorem 2.7 in~\cite{coates2}).
  Let $H=\Gal(F_\infty/F^{cyc})$ and $\Omega=\Gal(F_\infty/L^{cyc})$. The heart of the proof is the analysis of the following commutative diagram:-
  \begin{diagram}\label{extract1}
    0   &    \rTo & \SE^\Omega                  & \rTo & H^1(G_S(F_\infty),\Ep)^\Omega   & \rTo &  \bigoplus_{v\in S}J_v(F_\infty)^\Omega & \rTo & 0 \\
        &         &   \uTo_{f}                  &      &  \uTo_{g}                       &      &     \uTo_{h}                            &      &\\
    0   &    \rTo & \Sel(E/L^{cyc})             & \rTo & H^1(G_S(L^{cyc}),\Ep)           &\rTo  &    \bigoplus_{v\in S} J_v(L^{cyc})      & \rTo & 0\\
  \end{diagram}
  As shown in Propositions~\ref{globEC} and \ref{globEC*}, the map $g$ has finite kernel and cokernel. Note that for each prime $v$ of $F$, there are only finitely
  many primes of $L^{cyc}$ dividing it. By the inflation-restriction exact sequence, the kernel of $h$ is given by
  \begin{equation}
   \ker(h)=\bigoplus_{v\in S}H^1(\Omega_w,E)(p),
  \end{equation}
  where for each $v\in S$, $w$ is a fixed prime of $F_\infty$ above $v$. As before, we have
  \begin{align}\label{kernel(h)}
   H^1(\Omega_w,E)(p)\cong H^1(\Omega_w,\Ep)                    &\text{             if $w\nmid p$},\\
   H^1(\Omega_w,E)(p)\cong H^1(\Omega_w,\tilde{E}_{w,p^\infty}) &\text{             if $w\mid p$}.
  \end{align}
  Now $\Ep\cong\ZZ_p^2$ and $\tilde{E}_{w,p^\infty}\cong\ZZ_p$ for all $w\mid p$, so their Pontyragin duals are finitely generated $\ZZ_p$ -modules and in 
  particular finitely generated $\Lambda(\Omega_w)$-modules. As described in~\cite{howson1}, it follows that the cohomology groups in (\ref{kernel(h)}) have 
  finite $\ZZ_p$-corank. We now deduce from the snake lemma that $\ker(f)$ and
  $\coker(f)$ have finite $\ZZ_p$-coranks. By assumption, ${\cal C}(E/L^{cyc})$ is a
  finitely generated $\ZZ_p$-module, so $\Sel(E/F_\infty)^\Omega$ has finite $\ZZ_p$-corank. Since $\Gal(F_\infty/L)$ is pro-$p$, it follows that 
  $\Omega$ is also pro-$p$. Now ${\cal C}(E/F_\infty)$ is compact, so we may apply Nakayama's Lemma to deduce that ${\cal C}(E/F_\infty)$ is a 
  finitely generated $\Lambda(\Omega)$-module. But $\Omega$ is a subgroup of $H$, so ${\cal C}(E/F_\infty)$ is a finitely generated $\Lambda(H)$-module.
  Assume now that $\Sigma$ has no element of order $p$. 
  Let $H_0=\Gal(F_\infty/L^{cyc}(A_p))$, and recall that $R=\Gal(F_\infty/F(A_p))$. Note that $H_0$ has finite index in $H$, so as above we deduce that ${\cal
  C}(E/F_\infty)$ is finitely generated as a $\Lambda(H_0)$-module. Recall that ${\cal C}(E/F_\infty)$ is defined to be $\Lambda(\Sigma)$-torsion if it is
  torsion as a $\Lambda(R)$-module.
  Now $H_0$ has infinite index in $R$, so $\Lambda(R)$ is not finitely generated as a $\Lambda(H_0)$-module. Since $R$ has no element of order $p$, its Iwasawa
  algebra $\Lambda(R)$ contains no non-trivial zero-divisors. Hence we
  deduce that any finitely generated $\Lambda(R)$-module which is finitely generated as a $\Lambda(H)$-module must be $\Lambda(R)$- and hence
  $\Lambda(\Sigma)$-torsion, which proves the theorem.
 \end{proof} 

 In this section, we let $p=7$ and $F=\QQ(\mu_7)$. To illustrate the main result of the paper, we calculate the Euler characteristic of the Selmer group 
 of the elliptic curve
 \begin{equation}
  E:y^2+xy=x^3-x-1
 \end{equation}
 over several $7$-adic Lie extensions of $F$ the type which is considered in this paper. We first need some data about $E$ which we quote from
 ~\cite{coatessujatha1}: The curve $E$ has conductor $294$ and discriminant $\Delta=-2\cdot 3\cdot 7^2$; it has a point of order $7$ over $F$; it achieves good ordinary 
 reduction at the unique prime of $F$ 
 above $7$ and has split multiplicative reduction at the the two primes of $F$ above 2 and at the unique prime above $3$. It is shown in~\cite{fisher} that
 $\Sel(E/F)=0$ and that $\Sel(E/F^{cyc})=0$. It follows that ${\cal C}(E/F^{cyc})=0$, so in particular it is finitely generated as a $\ZZ_7$-module. Also, it is
 shown in~\cite{coatessujatha1} that
 \begin{equation}
  \rho_7(E/F)=1.
 \end{equation}
 Now let $E'$ be any elliptic curve defined over $F$ which has a point of order $7$ and which achieves good reduction at the unique prime of $F$ above $7$. 
 Let $A=E\times E'$, $F_\infty=F(A_{7^\infty})$ and
 $\Sigma=\Gal(F_\infty/F)$. Now $E$ has good reduction at all $v\nmid \{2,3,7\}$. The local Euler factor of $E$ at a prime $v\nmid 7$ of $F$ is
 \begin{equation}
     L_v(E,1)=
     \begin{cases}
      (1-2^{-3s})^{-1} & \text{if $v\mid 2$} \\
      (1-3^{-6s})^{-1} & \text{if $v\mid 3$} \\
      1         & \text{otherwise}.
     \end{cases}
 \end{equation}    
 Since ${\cal C}(E/F^{cyc})$ is finitely generated as a $\ZZ_7$-module, it follows from Theorem~\ref{finitegeneration} that ${\cal C}(E/F_\infty)$ is
 $\Lambda(\Sigma)$-torsion, so all the hypotheses of Theorem~\ref{EC} are satisfied. We deduce that
 \begin{equation}
  \chi(\Sigma,\Sel(E/F_\infty))=7^3.
 \end{equation}
 The situation is different when $E_{7^\infty}$ is not rational over the $7$-adic Lie extension $F_\infty$:- Let 
 \begin{equation}
  E':  y^2-xy+2y=x^3+2x^2,
 \end{equation}
 which has discriminant $\Delta '=-2^7 13$. Let $F_\infty=F(E'_{7^\infty})$ and $\Sigma=\Gal(F_\infty/F)$. Since $7\nmid\Delta '$, $E'$ has good reduction at
 the prime of $F$ above $7$, and it is easy to see that the point $(0,0)$ is of order $7$. We have ${\frak M}=\{2,13\}$. Using again above observations, it follows that the extension
 $F_\infty$ over $F$ satisfies the conditions of Theorem~\ref{EC}. We find that this time,
 \begin{equation}
  \chi(\Sigma,\Sel(E/F_\infty))=7^2.
 \end{equation}

\bibliography{references}

\end{document}